\RequirePackage{fix-cm}
\documentclass[smallextended]{svjour3}       
\smartqed  
\usepackage[utf8]{inputenc}
\usepackage[T1]{fontenc}
\usepackage[english]{babel}

\usepackage[usenames]{color}
\usepackage{colortbl}

\usepackage{amsmath}
\usepackage{multirow}

\usepackage{graphicx}

\DeclareMathOperator{\arcctg}{arcctg}
\begin{document}

\title{Tensor train optimization for mathematical model of social networks \thanks{This work is supported by the Russian Science Foundation (project No.~18-71-10044).}
}


\author{Kabanikhin S., \and Krivorotko O., \and Zhang~S., \and Kashtanova V. \and Wang~Y. 
}


\institute{Kabanikhin Sergey \at
              Institute of Computational Mathematics and Mathematical Geophysics SB RAS, Novosibirsk, Russia \\
              Novosibirsk State University, Novosibirsk, Russia \\
              Tel.: +7-913-910-34-78\\
              Fax: +7-383-330-87-83\\
              \email{ksi52@mail.ru}           
           \and
           Krivorotko Olga \at
              Institute of Computational Mathematics and Mathematical Geophysics SB RAS, Novosibirsk, Russia \\
              Novosibirsk State University, Novosibirsk, Russia \\
              Tel.: +7-983-303-30-83\\
              Fax: +7-383-330-60-46\\
              \email{krivorotko.olya@mail.ru}
           \and
       Zhang Shuhua \at
            Tianjin University of Finance and Economics, Tianjin, China\\
            \email{shuhua55@126.com}
       \and
       Kashtanova Victoriya \at
          Novosibirsk State University, Novosibirsk, Russia\\
          \email{vikakashtanova@ya.ru}
    \and
       Wang Yufang \at
          Tianjin University of Finance and Economics, Tianjin, China\\
          \email{wangyufangminshan@163.com}   
}

\date{Received: date / Accepted: date}

\maketitle

\begin{abstract}
The optimization algorithms for solving multi-parameter inverse problem for the mathematical model of parabolic equations arising in social networks, epidemiology and economy are investigated. The data fitting is formulated as optimization of least squares misfit function. Firstly, the tensor train decomposition approach is presented as global convergence algorithm. The idea of proposed method is to extract the tensor structure of the optimized functional and use it for optimization. Then the inverse problem solution is reached by implementation of the local gradient approach. The evident formula for the gradient of the misfit function is obtained. The inverse problem for the diffusive logistic mathematical model described online social networks is solved by combination of tensor train optimization and local gradient methods. The numerical results are presented and discussed.
\end{abstract}

\section{Introduction}
\label{intro}
The specifics of the dissemination of information in society and the development of socially significant diseases (tuberculosis, HIV/AIDS) depend on the region. However, the statistically calculated parameters (for example, the probability of the appearance of information in the social network, the rate of infection development, the mortality parameters, etc.) are average. The investigation of model populations, on the one hand, can lead to erroneous conclusions (since real populations are significantly heterogeneous) and, on the other hand, it cannot serve as a tool for assessing the current situation.

One of the most effective methods of monitoring and managing social and epidemiological processes is mathematical modeling, namely the development and identification of mathematical models that describe the processes of information dissemination in social networks and infections in the population. Such models are described by systems of differential equations, the coefficients of which characterize the distribution of information, population and disease development. To control information in social networks and epidemics in individual regions and economic processes, it is necessary to refine the model coefficients and initial data by some additional information (the inverse problem). One way to solve the problem of improving the coefficients and initial data is to reduce the inverse problem to a variational formulation, where the functional characterizes the quadratic deviation of the model data from the experimental ones. 

In this paper the inverse problem for the mathematical model of social networks based on parabolic equation is numerically investigated. The inverse problem is reduced to the optimization problem (Section~\ref{sec:3}). The minimization problem is solved by the tensor train decomposition and gradient methods (Section~\ref{sec:4}). The comparison and domain of applicability of investigated methods are analyzed and discussed (Section~\ref{sec:5}).

\section{Brief historical review}
\label{sec:2}
\textbf{Ordinary differential equations}. In recent years, online social networks, such as Twitter and Facebook, have become the main source of information exchange. A large amount of data available to researchers has increased interest in studying of the process of information  dissemination in online networks. One of the approaches to constructing mathematical models of social processes is the principles of constructing mathematical models of epidemiological processes (in particular, describing socially significant diseases, such as tuberculosis, HIV/AIDS) based on the chamber structure and probabilistic transitions between homogeneous groups~\cite{1} (mathematical models of tuberculosis~\cite{2},~\cite{3},~\cite{4}, HIV~\cite{5}). The spread of the virus in a homogeneous network based on the epidemiological model was investigated in the work~\cite{6}. Such mathematical models are based on systems of ordinary differential equations, the coefficients of which in many cases are unknown or given approximately. This leads to solving of inverse problems. In many models it is assumed that the social system is homogeneous, and individuals are in equal status (the same degree of spread, the probability of infection, etc.), which does not take into account the specifics of the process. This entails the development of more complex mathematical models.

\textbf{Partial differential equations.} To build a more complete picture of the development of social processes, it is necessary to take into account migration, age data and changes in time. Such models are described by partial differential equations (PDE). F.~Wang, H.~Wang and K.~Xu~\cite{8} proposed to use partial differential equations built on intuitive cyber-distance among online users to study both temporal and spatial patterns of information diffusion process in online social networks. A detailed review of mathematical models for social networks is given in~\cite{9}. The PDE-based models for online social networks in~\cite{10} are spatial dynamical systems that take into account the influence of the underly network structure as well as information contents to predict information diffusion over both temporal and spatial dimensions.
In paper~\cite{11} it is developed a non-autonomous diffusive logistic model with indefinite weight and the Robin boundary condition to describe information diffusion in online social networks. The model is validated with a real dataset from an online social network, Digg.com, and the simulation shows that the logistic model with the Robin boundary condition is able to more accurately predict the density of influenced users. 

\textbf{Stochastic differential equations.} Real social systems are always exposed to external influences that are not completely understandable or impossible for explicit models (enzymatic processes, energy needs, smoking, stress effects, information wars, etc.), and therefore there is a growing need to expand deterministic models to models, which cover more complex variations in dynamics. A natural continuation of the models of deterministic differential equations is a system of stochastic differential equations, where the corresponding parameters are modeled as suitable random processes, or stochastic processes are added to the equations of the motion system. An analysis of stochastic differential equations and numerical studies of the solution of a direct problem are given in the works of P.-L.~Lions~\cite{12},~\cite{13} and H.T.~Banks~\cite {14},~\cite{15}.

However, each social network has its own platform and structure, and therefore the parameters that characterize these indicators vary. For the best result of modeling and information control it is necessary to identify the coefficients and initial conditions (inverse problem) for each specific case. In papers of Kabanikhin S.I., Krivorotko O.I. and co-authors~\cite{16},~\cite{17} the inverse problem for ordinary differential equations describing the spread of tuberculosis and HIV in the regions of the Russian Federation and the inverse problem of immunology (intracellular dynamics of HIV)~\cite{18} were solved. An analysis of the structural identifiability of the mathematical models under study, which makes it possible to identify the region of the correctness of the inverse problem, is given in the paper~\cite{19}.

In the next Sections we concentrate on PDE-based mathematical models arises in social networks, epidemiology, immunology, economics.

\section{Inverse problem statement for PDE}
\label{sec:3}
The initial-boundary value problem for the PDE in parabolic type has the following form:
\begin{subequations}
\begin{align}
&\dfrac{\partial y_j}{\partial t} = d_j\dfrac{\partial^2 y_j}{\partial x^2} + \mu_j (y,\varphi), \quad & t\in(t_0,T), x\in(l,L);     \label{eqn: math_model}\\
&y_j(x,t_0)=\psi_j(x), \quad & x\in(l,L);   \label{eqn: init_data}  \\
&\dfrac{\partial y_j}{\partial x}\bigr|_{x=l}=\dfrac{\partial y_j}{\partial x}\bigr|_{x=L}=0,\quad & t\in(t_0,T)   \label{eqn: boundary_cond}
\end{align}
\end{subequations}
Here $y(x,t)=(y_1(t,x),\ldots,y_N(t,x))$, $\varphi(t)=(\varphi_1,...,\varphi_M)$ is a vector of unknown coefficients, $l>0$, $\mu_j$ is continuous function, $d_j>0$, $j=1,\ldots,N$.

Suppose that the additional information can be measured:
\begin{equation}\label{data}
y_j(x_i,t_k;q)=F^j_{ik}, \; k=1,\ldots,K, i=1,\ldots,N_x, j\subset J:=\{1,...,N\}.
\end{equation}

\textbf{The inverse problem~(\ref{eqn: math_model})-(\ref{eqn: boundary_cond}), (\ref{data})} consists in determination of coefficients and initial data $q=(\varphi(t), d, \psi(x))\in Q$ using additional information $F^j_{ik}$ (\ref{data}). Here $Q$ is the set of admissible solutions of the inverse problem.

The mathematical model~(\ref{eqn: math_model}) can describe the dynamic of social, economical and epidemiological processes. For example, in \textbf{social networks} $j=1$ and parameters and functions in problem~(\ref{eqn: math_model})-(\ref{eqn: boundary_cond}) have the following interpretations~\cite{11}:
\begin{itemize}
\item $y(x,t)$ is a density of influenced users with a distance of $x$ at time $t$;
\item $d$ represents the popularity of information which promotes the spread of the information through non-structure based activities such as search;
\item $\mu(y,\varphi) = r(t)y(x,t)\left(1-\frac{y(x,t)}{K_{cap}}\right)$ is a local growth function (death and birth);
\item $r(t)$ represents the intrinsic growth rate of influenced users with the same distance, and measures how fast the information spreads within the users with the same distance;
\item $K_{cap}$ is a carrying capacity, which is the maximum possible density of influenced users at a given distance;
\item $l$ and $L$ represent the lower and upper bounds of the distances between the source and other social network users;
\item $\psi(x)\geq0$ is the unknown initial density function. Each information has its own unique initial function.
\end{itemize}

The mathematical model (\ref{eqn: math_model}) arises in the modern \textbf{economy} and is called the spatial Solow model~\cite{23},~\cite{24} ($j=1$):
\begin{itemize}
\item $y(x,t)$ denotes the capital stock held by the representative household located at $x$ and time $t$;
\item $\mu(y,\varphi) = A(x,t)f(y(x,t)) - \varepsilon y(x,t)$;
\item $A(x,t)$ denotes the technological level at $x$ and time $t$;
\item $f(y(x,t))$ is the production function that assumed to be non-negative, increasing and concave;
\item $\varepsilon$ is the depreciation rate;
\item $\psi(x)\geq0$ is an initial capital distribution.
\end{itemize}

A PDE~(\ref{eqn: math_model}) for age-structured populations describes the \textbf{epidemiological} process and is called Kermack-McKendrick equation~\cite{22}:
\begin{itemize}
\item $y_j(x,t)\geq0$ represents the density of the population (infected and recovered individuals of different groups) of age $x$ and time $t$;
\item $\mu_j(y,\varphi)$ is a local growth function that describes the behavior of $j$-th group;
\item the sparsity coefficient $d_j$ characterizes the migration at $j$-th group;
\item $\psi_j(x)\geq0$ is an initial population at $j$-th group.
\end{itemize}

Inverse problem~(\ref{eqn: math_model})-(\ref{eqn: boundary_cond}), (\ref{data}) is ill-posed, and most often its solutions are unstable to small perturbations in the data~\cite{25}. One of the type of regularization of ill-posed problems is reducing inverse problems to a variational formulation. For this purpose we consider the misfit function
\begin{equation}\label{eq:min_problem}
J(q)=\dfrac{(T-t_0)\cdot (L-l)}{K\cdot N_x}\sum\limits_{j\subset J}\sum\limits_{k=1}^K\sum\limits_{i=1}^{N_x} |y_j(x_i,t_k;q)-F_{ik}^j|^2.
\end{equation}
that should to be minimized. Finally, the solution $q^\star$ of the inverse problem~(\ref{eqn: math_model})-(\ref{eqn: boundary_cond}), (\ref{data}) can be considered as the solution of an optimization problem~\cite{25}:
\[q^\star=\mbox{arg}\min\limits_{q\in Q} J(q).\]

\section{Optimization methods}\label{sec:4}
There are several methods for solving the minimization problem (\ref{eq:min_problem}), and they can be split into three groups:
local, global and hybrid optimization methods. An extensive review of the existing methods, including their classifications and properties can be found in~\cite{FEBS_parest}.
There are two groups of local optimization methods: the direct ones which do not require the derivative of the objective function
(e.g., Hooke-Jeaves method, Nelder-Mead method) and the gradient based methods (e.g., the Gauss-Newton- and Levenberg-Marquard methods).
Although the local optimization methods usually work fast and their convergence to a local minimum can be proved,
their main drawback is that they can miss the global minimum, especially in the case of a high dimensionality of the parameter space.
Global optimization methods can be used to explore rather large regions of the parameter space. However, they usually are slow and have no theoretical proof of convergence.
They can be formulated as stochastic (simulated annealing, evolutionary algorithms) or deterministic (covering methods) techniques~\cite{FEBS_parest} such as method of nonuniform coverings proposed by Yu.G.~Evtushenko for functions that comply with the Lipschitz condition~\cite{Evtushenko_1971}.
Evolutionary algorithms are usually most suitable for large search space because they escape local minima and are intrinsically parallel.
The covering methods require some prior information about the function and can locate the optima with the given accuracy. For example, the main idea of the method of nonuniform coverings consists in division of the solution set into subsets whose union coincides with the original set. The misfit function on different subsets has certain properties (satisfies the Lipschitz condition, convexity, existence of the second derivative, etc.), which speeds up the calculations~\cite{Evtushenko_2007,Evtushenko_2009}.
Hybrid methods are based on the following idea: global optimization is used to explore the parameter space to locate the starting points for further local optimization~\cite{17}.
To reduce the complexity of the parameter estimation task, there exist several support techniques such as constraining the parameter space, data smoothing, and others~\cite{Voit_PE_BN}. 

The majority of the state-of-the-art global optimization methods are based on tensor-representation of initial functional and applying of the tensor decomposition properties and tensor interpolations~\cite{ot-tt-2009,ot-ttcross-2010,osel-tt-2011}, that we use in numerical calculations and show its advantages.

\subsection{Tensor Train decomposition}

Optimization of the functional (\ref{eq:min_problem}) is a complex problem because it is a problem of global optimization for a multidimensional, usually, non-convex functional with large dimensional and huge number of local minimums. 

Tensor train (TT) approach~\cite{osel-tt-2011} based on decomposition of a multidimensional tensor $A\in R^{n_1\times n_2\times\ldots\times n_p}$ to the "train" of the \textit{carriages} tensors:
\begin{equation*}
\begin{aligned}
A(i_1;i_2;...;i_p)&=\sum\limits_{\alpha_0=1;...;\alpha_p=1}^{r_0;...;r_p}A_1(\alpha_0;i_1;\alpha_1)\cdot\ldots\cdot A_p(\alpha_{p-1};i_{p};\alpha_p). \\
\end{aligned}
\end{equation*}

Although this approach is metaheuristic and could not guarantee that the global optimum is reached, it uses structure of the functional and, generally, works faster and more robustly than other metaheuristic and stochastic methods~\cite{Zheltkova_Tyrtyshnikov_2018}.

\subsubsection{TT algorithm for optimization problem}
Consider the global minimization problem
\[q^\star = \arg \min_{q\in Q} J(q).\]
This problem could be transformed to an equivalent problem of the magnitude maximization of the continuous and monotonous mapping of $J$ to the interval $[0;+\infty)$:
\begin{equation}\label{minimization_pr}
q^\star = \arg \min_{q\in Q} |g(q)|,\quad g(q) = \arcctg \{J(q)\}.
\end{equation}

Here $Q$ is a initial domain of the function $J$ and can be represented as a $p$-dimensional parallelepiped ($p=M+N+1$). So $q=(q_1,\ldots,q_p)^T$ is a vector and, suppose, the each element $q_j$ lies in the interval $[a_j,b_j]$, $j=1,\ldots,p$. Introduce the uniform grid on each interval $[a_j,b_j]$ with fix step $h_j = (b_j-a_j)/(n-1)$,  $j=1,\ldots,p$ and $n$ nodes in each direction. 

Misfit function values on the grid form tensor $G\in R^{n\times\ldots\times n}$ with elements:
\[G(i_1,\ldots,i_p) = g(q_1^{(i_1)},\ldots, q_p^{(i_p)}).\]
Here $q_j^{(i_j)}$ is a $i_j$-th point on the grid for element $q_j$. Introduce vector of indexes $i=(i_1,\ldots,i_p)$ and rewrite the formula as
\[G(i) = g(q^{i}).\]
Then discrete minimization problem is represented in the follows form:
\begin{equation}\label{TT_min_pr}
i=\arg \min_{j=1,\ldots,n^p} |G(j)|.
\end{equation}
If the grid is fine enough then the solutions of (\ref{minimization_pr}) and (\ref{TT_min_pr}) are expected 
to be close.

The problem (\ref{TT_min_pr}) consists in finding the maximal in magnitude element of a $p$-dimensional tensor $ G \in R^{n \times \dots \times n}$. For example, if a number of coefficients is equal to 20 and a number of mesh nodes in space $x$ is equal to 100, then the number of unknown parameters is equal to 121. If $n=100$, then the number of elements of tensor $G$ is equal to $100^{121}$.

In this case, to reduce the problem complexity,  the technique based on the TT-cross interpolation machinery, which exploits the matrix cross interpolation algorithm~\cite{gt-psa-1995,gtz-psa-1997,tee-cross-2000} applied to heuristically selected submatrices in the unfolding matrices of the given tensor, could be used. The obtained method, named TT algorithm, takes only $O(pnr_{max}^3)$ arithmetic operations, $O(pnr_{max}^2)$ function calculations and $O(nr_{max})$ local optimizations, where $r_{max}$ is the maximum rank of used tensors.





The TT global optimization method iteratively performs the following steps:
\begin{itemize}
 \item already inspected points are used to generate submatrices of the unfolding matrices;
 \item these submatrices are approximated by the matrix cross approximation method with rank bounded from above by $r_{max}$;
 \item the interpolation points and local minimums in their vicinity (projected to the grid) are used to form 
new sets of ``hopefully better'' points;
 \item the sets of points are extended by the points from ``neighboring'' unfolding matrices 
and by $r_{max}$ points considered as the best of all inspected values.
\end{itemize}

\begin{figure}[!ht] 
\begin{center}
\includegraphics[width=0.9\textwidth]{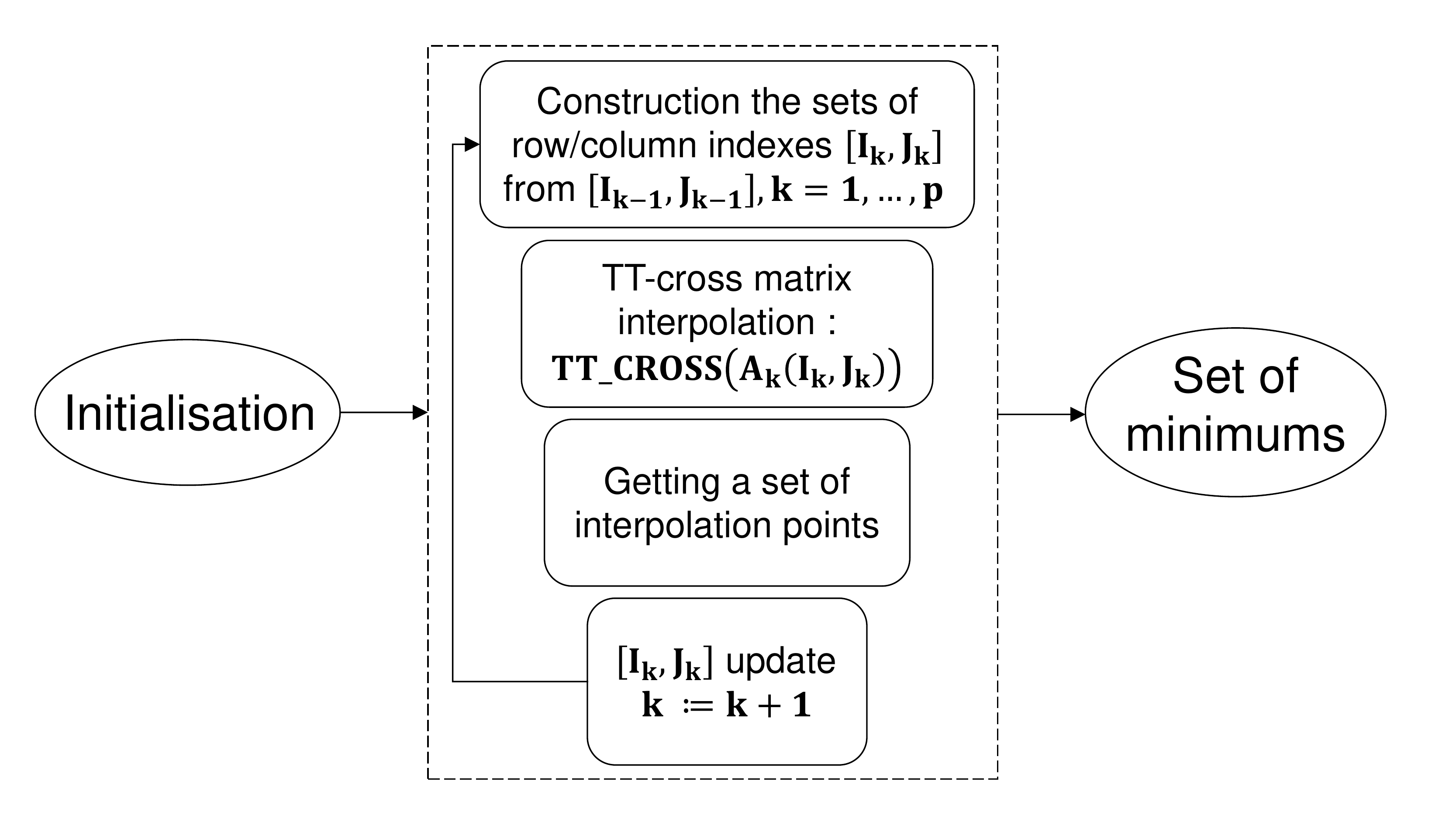}
\caption{The general scheme of TT decomposition algorithm.}
\label{fig:TT-sheme}      
\end{center}
\end{figure}

The general scheme of this method is presented by Fig.~\ref{fig:TT-sheme}, more detailed description of algorithm is in~\cite{Zheltkova_Tyrtyshnikov_2018,Zheltkov_Oferkin_2013}.

\subsection{Gradient method}
Most parts of gradient methods are reduced to iteration process of the following form:
\[ q_{n+1} = q_n - \alpha_n J^\prime(q_n), \]
where $\alpha_n$~is a descent parameter, $J^\prime(q_n)$~is a gradient of misfit function at point $q_n$. Note that convergence of gradient method depends on initial approximation $q_0$ and as a result has the local convergence theorem~\cite{25}.

The dependence on vector of parameters $q$ notes at $y(x,t;q):=y(x,t)$. Function $\hat{y}: (0,T) \mapsto H^2(l,L)$ is a mapping, associated by the function $y(x,t)$, that is, $[\hat{y}(t)](x) := y(x, t)$. Suppose that $\psi$, $\psi+\delta\psi\in H^2(l,L)$, $\delta y(x,t;\delta q):=y(x,t;q+\delta q) - y(x,t;q)$, where $y(x,t;q)\in C(0,T;H^2(l,L))$ is a solution of problem~(\ref{eqn: math_model})-(\ref{eqn: boundary_cond}) with $\psi\in H^2(l,L)$. Then the deviation $\delta y:=\delta y(x,t;\delta q)$ satisfies the following initial boundary-value problem with an accuracy up to the terms of
order $o(|\delta q|^2)$:
\begin{eqnarray}\label{eq:perturbation_y}
\left\{\begin{array}{ll}
L\delta y_j=\dfrac{\partial \delta y_j}{\partial t} - d_j\dfrac{\partial^2 \delta y_j}{\partial x^2} - a_j(\delta y,\delta\varphi,y,\delta d)=0, \quad & t\in(0,T), x\in(l,L);\\
\delta y_j (x,0) = \delta \psi_j(x), \quad & x\in(l,L);\\
\left.\frac{\partial \delta y_j}{\partial x}\right|_{x=l}=\left.\frac{\partial \delta y_j}{\partial x}\right|_{x=L}=0,\quad & t\in(0,T).
\end{array}\right.
\end{eqnarray}
Here $a_j(\delta y,\delta\varphi,y,\delta d) = \sum\limits_{n=1}^{N} P_{jn} \delta y_n + \sum\limits_{m=1}^{M} R_{jm} \delta\varphi_m + \delta d_j\dfrac{\partial^2 y_j}{\partial x^2}$, $j=1,\ldots N$, $P=\left\{\dfrac{\partial \mu (y,\varphi)}{\partial y}\right\}$ and $R=\left\{\dfrac{\partial \mu (y,\varphi)}{\partial \varphi}\right\}$ are Jacobi matrices of vector function $\mu$ at spaces $\mathbf{R}^N\times \mathbf{R}^N$ and $\mathbf{R}^N\times \mathbf{R}^M$, respectively.

The following proposition is hold:

\textbf{\underline{Proposition 1.}} The gradient of misfit function~(\ref{eq:min_problem}) has the following form:
\begin{equation}\label{grad_func}
J'(q)=\Bigl(-\int\limits_{l}^L R^T\Psi(x,t)\,dx, -\int\limits_{0}^L\int\limits_0^T \dfrac{\partial^2 y}{\partial x^2}(x,t)\Psi(x,t)\,dt dx, -\Psi(x,0)\Bigr)^T,
\end{equation}
where vector-function $\Psi(x,t)$ satisfies to the adjoint problem ($j=1,\ldots,N$):
\begin{equation}\label{eq:adj_pr}
\left \{ 
\begin{aligned}
&\dfrac{\partial\Psi_j}{\partial t} = - d_j \dfrac{\partial^2\Psi_j}{\partial x^2} - \sum\limits_{n=1}^{N} P_{nj}\Psi_n + b_j(x,t), \quad &t\in(0,T), x\in(l,L);  \\
&\Psi_j(x,T)=0, \quad &x\in(l,L);  \\
&\left.\frac{\partial \Psi_j}{\partial x}\right|_{x=l}=\left.\frac{\partial \Psi_j}{\partial x}\right|_{x=L}=0,\quad &t\in(0,T).
\end{aligned}\right.
\end{equation}
Here $b_j(x,t) = \sum\limits_{k=1}^K \sum\limits_{i=1}^{N_x} \int\limits_0^T\int\limits_l^L 2\left(y_j(x,t;q)-F^j_{ik}\right)\delta(t-t_k)\delta(x-x_i)\,dxdt$, $\delta(t-t_k)$ is a Dirac delta-function.

\textbf{\underline{Proof.}} Consider the misfit function~(\ref{eq:min_problem}) variation and apply the rule $a^2-b^2=(a-b)(a+b)$:
\begin{equation}\label{deviationJ}
\left.\begin{aligned}
\delta J&:=J(q+\delta q) - J(q) = \\
&=\sum\limits_{k=1}^K \sum\limits_{i=1}^{N_x} 2\left(y(x_i,t_k;q)-F_{ik}\right)\delta y(x_i,t_k;q) + \sum\limits_{k=1}^K \sum\limits_{i=1}^{N_x} |\delta y(x_i,t_k;\delta q)|^2=\\
&=\sum\limits_{k=1}^K \sum\limits_{i=1}^{N_x} \int\limits_0^T\int\limits_l^L 2\left(y(x,t;q)-F_{ik}\right)\delta y(x,t;q)\delta(t-t_k)\delta(x-x_i)\,dxdt +\\
&+ \sum\limits_{k=1}^K \sum\limits_{i=1}^{N_x} |\delta y(x_i,t_k;\delta q)|^2.
\end{aligned}\right.
\end{equation}
Consider the scalar product in $L_2((0,T)\cup (l,L))$ space:
\begin{equation}\label{eq:scalar}
\left< L\delta y_j,\Psi \right> = \left<\dfrac{\partial \delta y_j}{\partial t},\Psi_j\right> - d_j\left<\dfrac{\partial^2 \delta y_j}{\partial x^2},\Psi_j\right> - \left<a_j(\delta y,\delta\varphi,y,\delta d),\Psi_j\right>.
\end{equation}
Write in details each term of the equation~(\ref{eq:scalar}). Using differentiation by parts and initial conditions of the direct~(\ref{eqn: init_data}) and adjoint~(\ref{eq:adj_pr}) problem we get the following expression for the first term:
\[\left<\dfrac{\partial \delta y_j}{\partial t},\Psi_j\right> = \int\limits_{0}^L\int\limits_0^T \dfrac{\partial \delta y_j}{\partial t}\Psi_j \,dtdx = -\left<\delta y_j,\dfrac{\partial \Psi_j}{\partial t}\right> - \int\limits_l^L \delta\psi_j(x)\Psi_j(x,0)\,dx.
\]
The second term in~(\ref{eq:scalar}) after twice usage of differentiation by parts and boundary conditions of the direct~(\ref{eqn: boundary_cond}) and adjoint~(\ref{eq:adj_pr}) problems has the form:
\[d_j\left<\dfrac{\partial^2 \delta y_j}{\partial x^2},\Psi_j\right> = -d_j\int\limits_{0}^L\int\limits_0^T \dfrac{\partial^2 \delta y_j}{\partial x^2}\Psi_j\, dt dx = -\left<\delta y_j,d_j \dfrac{\partial^2 \Psi_j}{\partial x^2}\right>.
\]
The last term in~(\ref{eq:scalar}) rewrite as follows using the linearity of a scalar product:
\begin{equation*}
\begin{aligned}
	&\left<a_j(\delta y,\delta\varphi,y,\delta d),\Psi_j\right> = \left<\sum\limits_{n=1}^{N} P_{jn} \delta y_n,\Psi_j\right> + \left<\sum\limits_{m=1}^{M} R_{jm} \delta\varphi_m,\Psi_j\right> +\\
	&+ \left<\delta d_j\dfrac{\partial^2 y_j}{\partial x^2},\Psi_j\right> = \left<\delta y_j,\sum\limits_{n=1}^{N} P_{nj}\Psi_n\right> + \left<\delta\varphi_j,\sum\limits_{m=1}^{M} R_{mj} \Psi_m\right> + \delta d_j\left<\dfrac{\partial^2 y_j}{\partial x^2},\Psi_j\right>.
\end{aligned}
\end{equation*}

Collect all equations to expression~(\ref{eq:scalar}), use the adjoint problem~(\ref{eq:adj_pr}) and note that formula~(\ref{eq:scalar}) is equal to zero we get the following equation for $j=1\ldots,N$:
\begin{equation*}
\begin{aligned}
	&b_j(x,t)\delta y_j(x,t;q) = - \int\limits_l^L \delta\psi_j(x)\Psi_j(x,0)\,dx -\\
	& -\delta d_j\int\limits_0^T\int\limits_l^L\dfrac{\partial^2 y_j}{\partial x^2}\Psi_j\,dxdt - \int\limits_0^T \delta\varphi_j\int\limits_l^L\sum\limits_{m=1}^{M} R_{mj} \Psi_m\, dxdt.
\end{aligned}
\end{equation*}
Using formula~(\ref{deviationJ}) note that:
\begin{equation*}
	\sum\limits_{j=1}^N b_j(x,t)\delta y_j(x,t;q) = \delta J - \sum\limits_{k=1}^K \sum\limits_{i=1}^{N_x} |\delta y(x_i,t_k;\delta q)|^2.
\end{equation*}

Suppose that $\sum\limits_{k=1}^K \sum\limits_{i=1}^{N_x} |\delta y(x_i,t_k;\delta q)|^2 \approx o(\Vert \delta q\Vert^2)$. Then after the comparison with the Freshet derivative formula for misfit function
\[\delta J = \left< J^\prime, \delta q\right> + o(\Vert \delta q\Vert^2),\]
we get the gradient $J^\prime(q)$~(\ref{grad_func}) from space $\mathbf{R}^{M+N+1}$.

\subsubsection{Algorithm of the gradient method of minimum errors}
The type of the gradient method depends on descent parameter $\alpha_n$. Propose the algorithm of a gradient method of minimum errors~\cite{17} as follows:
\begin{enumerate}
    \item Set an initial approximation vector $q_0$ and stopping parameter $\varepsilon>0$. Suppose that we have $q_n$. Show how to get the next approximation $q_{n+1}$.
    \item Check the stop condition: if $J(q_{n}) < \varepsilon$, then $q_{n}$ is an approximate solution of the inverse problem~(\ref{eqn: math_model})-(\ref{eqn: boundary_cond}), (\ref{data}). Otherwise, go to step 3.
    \item Solve the direct problem~(\ref{eqn: math_model})-(\ref{eqn: boundary_cond}) for a given set of the parameters $q_n$ by an explicit finite difference scheme of second order approximation and get $y_j(x_i,t_k;q_n)$, $i=1,\ldots,N_x$, $k=1,\ldots,K$, $j\subset J$.
    \item Solve the adjoint problem~(\ref{eq:adj_pr}) by an explicit finite difference scheme of second order approximation and get the solution $\Psi_j(x,t)$, $j\subset J$.
    \item Determine the gradient of misfit function $J(q_n)$ by the formula~(\ref{grad_func}).
    \item Calculate the descent parameter $\alpha_n = 2J(q_n)/\Vert J^\prime(q_n) \Vert$ for the minimum errors gradient method.
    \item Calculate the next approximation $q_{n+1} = q_n - \alpha_n J^\prime(q_n)$ and go to step 2.
\end{enumerate}

\section{Numerical solution of inverse problem for the mathematical models of social network}\label{sec:5}
Apply the proposed optimization algorithms to the inverse problem for the mathematical diffusive logistic model of type~(\ref{eqn: math_model}) arises in online social networks~\cite{8} (number of equations $N=1$):
\begin{eqnarray}\label{eqh:social_network}
\left\{\begin{array}{ll}
\dfrac{\partial y}{\partial t} = d\dfrac{\partial^2 y}{\partial x^2} + r(t)y(x,t)\left(1-\dfrac{y(x,t)}{K_{cap}}\right), \quad & t\in(1,T), x\in(l,L);\\
y(x,1)=\psi(x), \quad & x\in(l,L);  \\
\dfrac{\partial y}{\partial x}\bigr|_{x=l}=\dfrac{\partial y}{\partial x}\bigr|_{x=L}=0,\quad & t\in(1,T).
\end{array}\right.
\end{eqnarray}
\begin{eqnarray}\label{eqh:IP_social_network}
y(x_i,t_k) = F_{ik},\quad i=1,\ldots N_x, k=1,\ldots,K.
\end{eqnarray}
The description of all parameters and functions is presented in Section~\ref{sec:3}. Here we put
\[r(t) = \dfrac{\beta_2}{\beta_1} - e^{-\beta_1(t-1)}\left(\dfrac{\beta_2}{\beta_1} - \beta_3\right).\]
Since the initial density of influenced users $\psi(x)$ is unique for fixed social network and information type we get $\psi_i$, $i=1,\ldots,N_x$, from paper~\cite{8} (describe the situation of information network Digg.com) as an ''exact'' initial condition and then approximate it on $(l,L)$.
\begin{table}
\caption{The description and values of parameters of the mathematical model~(\ref{eqh:social_network}).}
\label{TAB:parameter_info}
\begin{tabular}{lll}
\hline\noalign{\smallskip}
Symbol & Description & Average value \\
\noalign{\smallskip}\hline\noalign{\smallskip}
$d$ & The popularity of information & 0.01 \\
$K_{cap}$ & A carrying capacity & 25 \\
$\beta_1$ & The rate of decline in information over time & 1.5 \\
$\beta_2$ & The residual speed & 0.375 \\
$\beta_3$ & Initial growth rate of the number of influenced users & 1.65 \\
\noalign{\smallskip}\hline
\end{tabular}
\end{table}

The inverse problem~(\ref{eqh:social_network})-(\ref{eqh:IP_social_network}) consists in determination of vector-function $q=(d, K_{cap}, \beta_1, \beta_2, \beta_3, \psi_1,\ldots,\psi_{N_x})$ from~(\ref{eqh:social_network}) using additional information~(\ref{eqh:IP_social_network}). For getting the synthetic inverse problem data~(\ref{eqh:IP_social_network}) we set an ''exact'' solution $q_{ex}$ according with table~\ref{TAB:parameter_info} and solve the direct model~(\ref{eqh:social_network}) with initial condition like in~\cite{8} by an explicit finite difference scheme of the second order approximation. We set the distance from source of information $N_x=6$ on friendship interval $l=1$, $L=6$ and measured the density of influenced users every hour from $t_1=5$ to $t_6=10$, $K=6$, during $T=24$ hours.

To obtain the best results we combine the algorithms according the pipeline: a global optimization algorithm to get the good initial approximation, then a local optimization method to get a result.

\begin{figure}[!ht] 
\begin{center}
\includegraphics[width=0.8\textwidth]{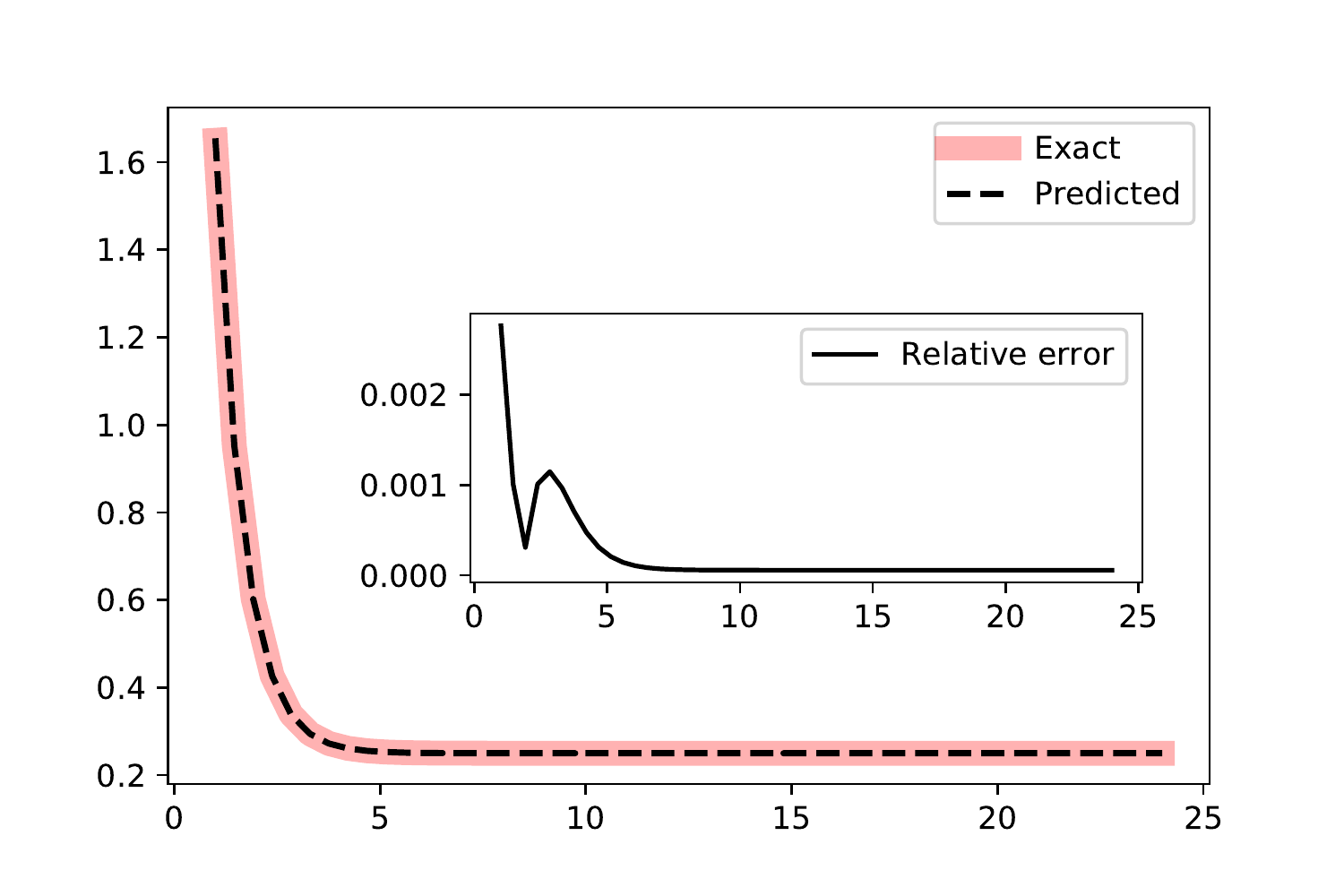}
\caption{The exact and reconstruction functions $r(t)$ with parameters from the table~\ref{TAB:parameter_reconstr} and the relative error $E(r)$.}
\label{fig:r(x)}      
\end{center}
\end{figure}

\begin{figure}[!ht] 
\begin{center}
\includegraphics[width=0.8\textwidth]{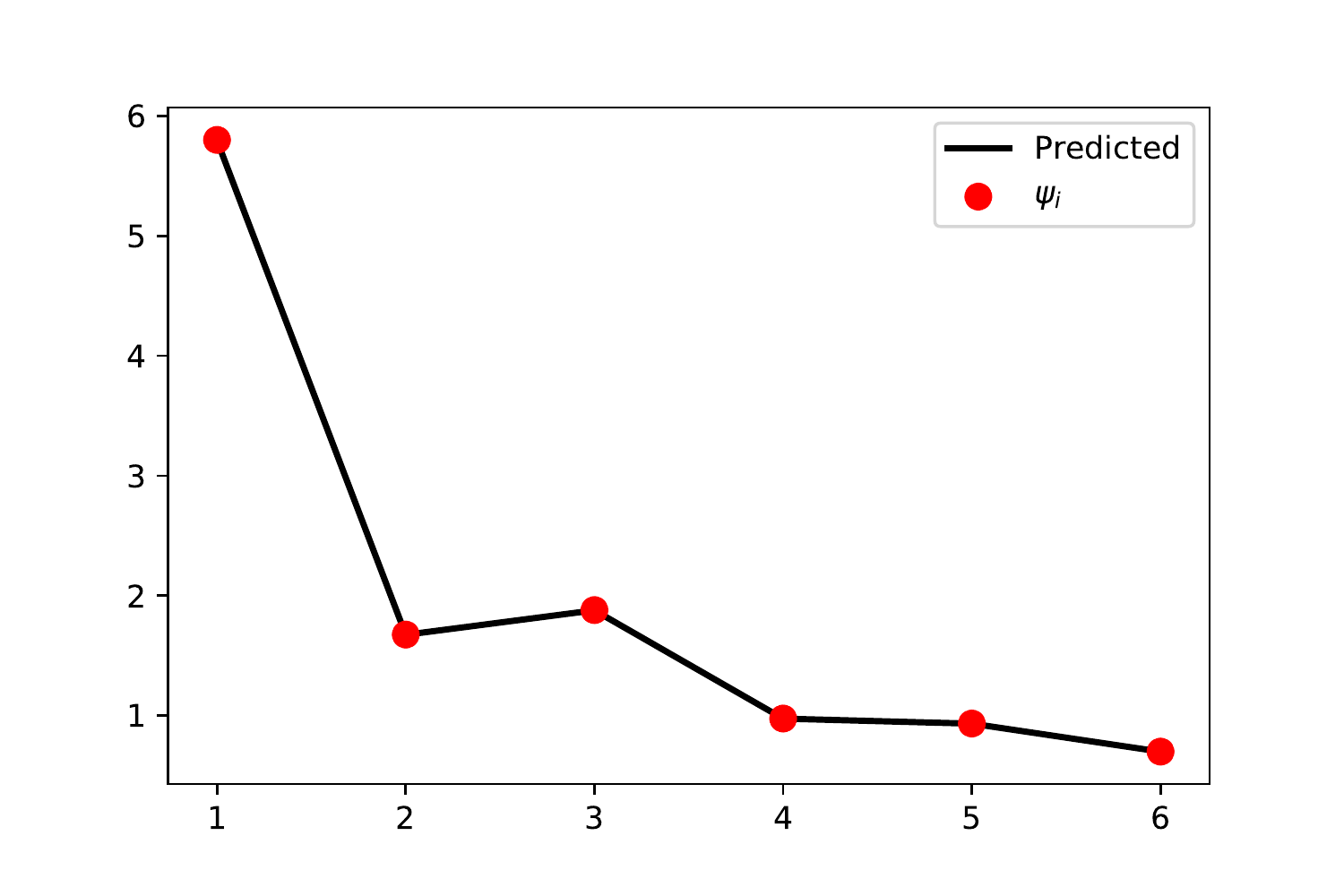}
\caption{Reconstruction of function $\psi(x)$ and ''exact'' points $\psi_i$, $i=1\ldots,6$. Here the relative error $E(\psi) = 9.5\cdot 10^{-4}$.}
\label{fig:phi(x)}      
\end{center}
\end{figure}

Solving the inverse problem under this conditions by the combined method, we got the following reconstructions of functions $r(t)$ and $\psi(x)$ (fig.~\ref{fig:r(x)}-\ref{fig:phi(x)}). Introduce the relative error for reconstruction parameters
\[E(r) = \frac{\Vert r_{pred} - r_{ex}\Vert_{L_2}}{\Vert r_{ex}\Vert_{L_2}}.\]

\begin{figure}[!ht] 
\begin{center}
\includegraphics[width=0.8\textwidth]{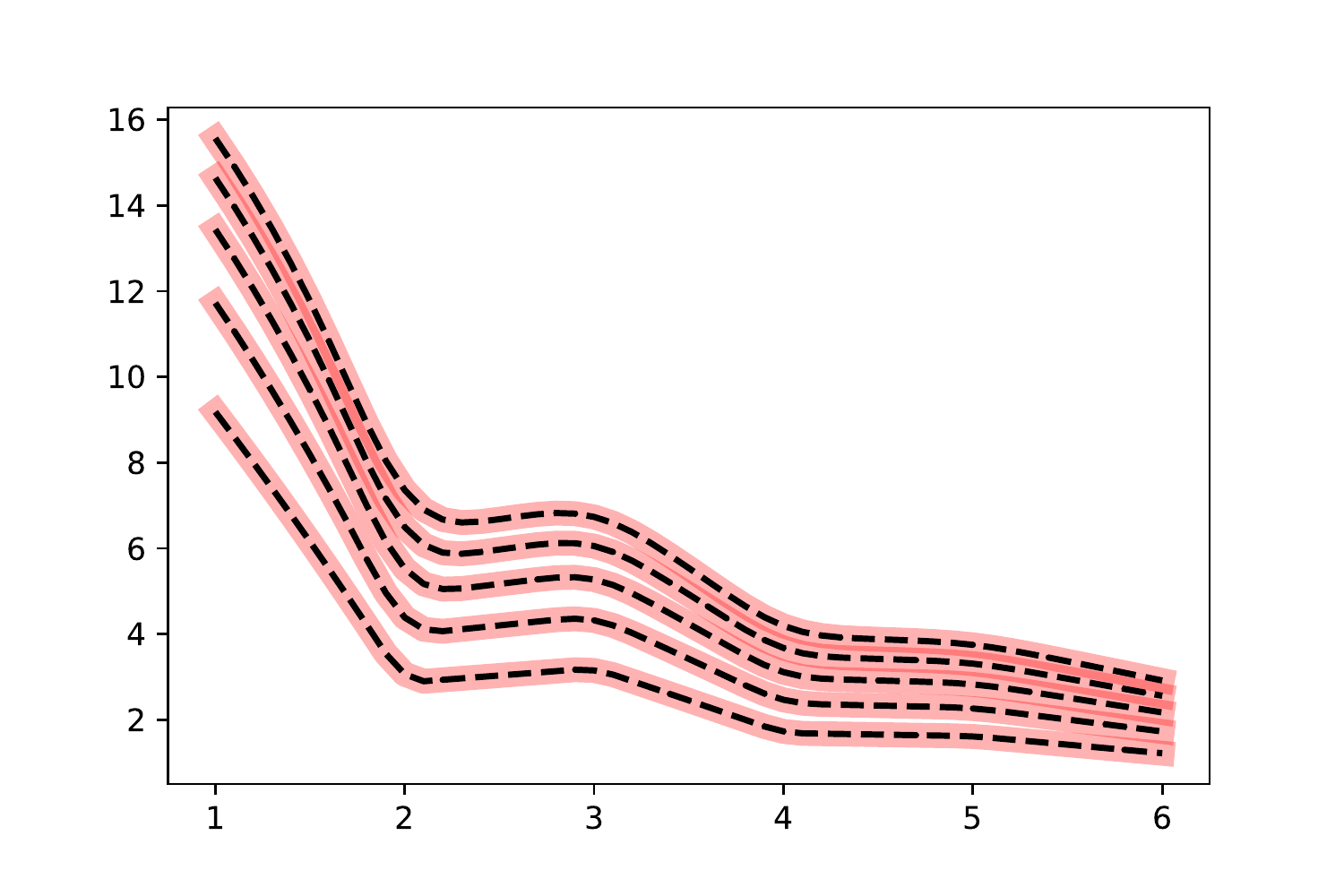}
\caption{Reconstruction of the direct problem for model~(\ref{eqh:social_network}) with parameters from the table~\ref{TAB:parameter_reconstr}. Solid width lines denote the real observations for the density of influenced users for a variety of time periods, dashed lines~-- the predicted density of influenced users.}
\label{fig:DP}      
\end{center}
\end{figure}

Figure~\ref{fig:DP} illustrates the predicting results for an
example news story with the proposed model, where the $x$-axis is the distance between users, while the $y$-axis represents the density of influenced users within each distance.
The solid width lines denote the actual observations for the density of influenced users for a variety of time periods (i.e., 1-hour, 2-hours, 3-hours, 4-hours and 5-hours), while the dashed lines illustrate the predicted density of influenced users by the model for reconstructed parameters $q$. As we can see, the proposed model is able to accurately predict the density of influenced users with different distance over time (see the real measurements in~\cite{8}).

Table~\ref{TAB:parameter_reconstr} gives the numerical value of parameters reconstruction during the numerical calculations with its relative approximation errors $E(q)$.

\begin{table}
\caption{The results of inverse problem~(\ref{eqh:social_network})-(\ref{eqh:IP_social_network}).}
\label{TAB:parameter_reconstr}
\begin{tabular}{llllll}
\hline\noalign{\smallskip}
\multirow{2}{*}{Symbol} & \multirow{2}{*}{Exact value} & \multicolumn{2}{c}{TT method} &  \multicolumn{2}{c}{Gradient method }\\
& & Value & Relative error & Value & Relative error \\
\noalign{\smallskip}\hline\noalign{\smallskip}
$d$ & 0.01 & 0 & 1 & 0.01 & $2.93\cdot 10^{-4}$ \\ 
$K_{cap}$ & 25 & 24.5 & 0.02 & 25 & $1.7\cdot 10^{-5}$ \\ 
$\beta_1$ & 1.5 & 1.73 & 0.154 & 1.504 & $2.7\cdot 10^{-3}$ \\ 
$\beta_2$ & 0.375 & 0.449 & 0.197 & 0.376 & $2.66\cdot 10^{-3}$ \\ 
$\beta_3$ & 1.65 & 1.73 & 0.05 & 1.654 & $2.76\cdot 10^{-3}$ \\ 
\noalign{\smallskip}\hline
\end{tabular}
\end{table}

\section{Conclusion}
The combined optimization algorithm for solving of multi-parameter inverse problem for the mathematical model of PDE of parabolic type arising in social networks, epidemiology and economy is investigated. The inverse problem consists in identification of coefficients in PDE and an initial condition of the initial boundary value problem for PDE using additional measurements of the solution of the direct problem in fixed points of one-dimensional space and time.
Considered inverse problem is ill-posed, i.e. the solution of inverse problem is non-unique and is unstable. We reduce the inverse problem to minimization of the least squares misfit function. There exists a wide class of optimization methods for multi-parameter minimization problem. We choose the tensor train decomposition approach. The idea of proposed method is to extract the tensor structure of the optimized functional and use it for optimization.
For more accurate reconstruction we apply the local gradient method of minimum errors. The approximate gradient of the misfit function includes the solution of an adjoint problem.

For numerical experiments the inverse problem for the diffusive logistic mathematical model described online social networks is solved by combination of tensor train optimization and local gradient methods.
All calculations were implemented on Python on the Google Cloud Platform with Virtual Machine with 15 Gb of RAM. It took about 18 million of function evaluations for tensor train method and about 3 thousand of function evaluations for gradient method. However, taking into account that in the tensor train method the functions are calculated in parallel the calculations took less machine time than in the gradient method.

Nevertheless we reached a good result by using the combined method, we could improve the quality of parameters reconstruction with the tensor train method by using a more detailed grid ($n = 100$). But, in this case, we need to use a more powerful computing machines.




\end{document}